\documentclass{elsarticle}
\usepackage{etoolbox}
\makeatletter
\patchcmd{\ps@pprintTitle}
  {Preprint submitted}
  {Published}
  {}{}
\makeatother

\def\ee{\end{eqnarray*}}
\def\be{\begin{eqnarray*}}
\def\bee{\end{eqnarray}}
\def\bbe{\begin{eqnarray}}

\newtheorem{Theorem}{Theorem}[section]
\newtheorem{Lemma}{Lemma}[section]

\def\e{\mathrm e}
\def\p{\partial }

\def\x{\mathbf x}
\def\y{\mathbf y}
\def\u{\mathbf u}
\def\w{\mathbf w}
\date{}
\journal{JMAA 431(2015) 1-21}
\begin{document}
\title{Steady-state bifurcation analysis of  a strong nonlinear  atmospheric  vorticity   equation
}

\author{Zhi-Min Chen
}

\address{School of Mathematics and Computational Science, Shenzhen University, Shenzhen, 
China\\
 Ship Science, University of Southampton, Southampton SO17 1BJ, UK}

\begin{abstract}
 The atmospheric  vorticity   equation  studied in the present paper is a simplified form of  the atmospheric circulation model introduced by Charney and DeVore [J. Atmos. Sci. 36(1979), 1205--1216] on the existence of multiple steady states to the understanding of the persistence of atmospheric blocking.  The fluid motion defined by the equation  is driven by a zonal  thermal forcing and an Ekman  friction forcing  measured by $\kappa$.
It is proved  that the   steady-state solution is globally unique for  large $\kappa$ values  while  multiple  steady-state solutions branch off  the basic steady-state solution  for $\kappa<\kappa_{\rm{crit}}$ where the critical value  $\kappa_{\rm{crit}}$ is less than one.
 Without involvement of viscosity, the equation has strong non-linear property   as its non-linear  part contains  the highest order derivative term.  Steady-state bifurcation analysis is essentially based on the compactness,   which can be simply  obtained for  semilinear  equations such as the Navier--Stokes equations but is not available for the strong nonlinear  vorticity  equation in the Euler formulation. Therefore the Lagrangian formulation of the equation  is employed  to gain the required compactness.

\

\noindent\textbf{Keywords:} atmospheric  vorticity   equation, steady-state bifurcation, Lagrange formulation, strong non-linear  equation

\

\noindent\textbf{Mathematics Subject Classification:} 35B32, 35B35, 35Q35, 86A10, 76B03

\end{abstract}
\maketitle

\section{Introduction} In an effort to describe the mechanism of   atmospheric  blocking phenomena,
 Charney and DeVore \cite{CD} introduced a two-dimensional quasi-geostrophic vorticity   equation and used a three mode truncation model to show heuristically the  existence of   multiple steady-state  solutions due to  non-linear  interaction of  zonal thermal forcing, Ekman layer energy dissipation  and  topography wave. Amongst them, a stable steady  state with weak zonal disturbance  describes the blocking phenomena. The numerical simulations of
the multiple  steady-state solutions of the  quasi-geostrophic vorticity   equations  originated from  Charney and DeVore \cite{CD} and have been extensively studied (see, for example,  
Eert \cite{HE}, Ierley  and  Sheremet \cite{IS}, Jiang et al. \cite{Jiangal}, Legras and   Ghil \cite{LG}, Pedlosky \cite{P},   Pierrehumbert  and  P. Malguzzi \cite{PM}, Primeau \cite{PF},  Rambaldi and   Mo \cite{RM}, Tung  and  Rosenthal \cite{TR}, Holloway  and   Yoden \cite{Y,Yo}) in the area of atmospheric science. However, the rigorous analysis supporting the multiple steady-state phenomenon is still lacking.

In the present paper, we are interested in the following atmospheric  vorticity   equation simplified from  Charney and DeVore \cite{C,CD}
 \bbe \label{nb1}
\frac{\partial \Delta \psi}{\partial t} +(\nabla \times \psi)\cdot \nabla (\Delta \psi)=-
\kappa \Delta (\psi-\psi^*)
\bee
with a flat  topography and the absence of the Coriolis force.  Here $\nabla$ is the gradient operator, $\Delta$ is the Laplacian,    $\psi$ is an unknown stream function,  $\kappa$ is an Ekman dissipative number, $\kappa \Delta \psi^*$ is an external thermal forcing and the vortex $\nabla \times \psi = ( -\partial_{x_2} \psi,\, \partial_{x_1} \psi)$.

This  is a strong nonlinear  third-order partial differential equation. If $\omega$ is employed to represent the vorticity $\Delta\psi$, the equation (\ref{nb1}) can be  rewritten in  the non-local form
 \be \label{nb1}
\frac{\partial \omega}{\partial t} +(\nabla \times \psi )\cdot \nabla \omega =-
\kappa (\omega-\omega^*)
\ee
due to involvement of the integral equation $\psi=\Delta^{-1}\omega$.  For the existence and singularities of evolutionary solutions to   related non-local equations, one may consult C\'ordoba et al. \cite{Cor} and Dong \cite{Dong}.

When $\kappa=0$, the equation (\ref{nb1}) reduces to the Euler equation. Thus the equation (\ref{nb1}) is the Euler equation with dissipation (see, for example,  \cite{Il1993}). The existence of a steady  state and the uniqueness of small steady  state  for the equation (\ref{nb1}) were  obtained by Wolansky \cite{Wol1} and Ilyin \cite{Il1993}. A more general form of  the equation (\ref{nb1}) is  known as the Stommel--Charney model \cite{Titi1988,Charney1955,Hauk1997,Stommel1948}, when the fluid motion involves the Coriolis force represented  the $\beta$ plane approximation in  middle latitudes. The existence of a steady  state and the uniqueness of small steady  state  for the Stommel--Charney model were obtained by   Barcilon et al.  \cite{Titi1988} and Hauk \cite{Hauk1997}.

However, the uniqueness may no long valid for large forcing  and  multiple   steady  states may coexist.  The purpose of present paper is to show  the existence of multiple steady-state solutions   of (\ref{nb1}) with respect to a parameter range of $\kappa$ and  the  zonal thermal forcing
\bbe\label{n3} \kappa \Delta \psi^*= -\kappa \cos x_2\,\,\mbox{ with } \psi^*=\cos x_2,
\bee
 employed  in \cite{CD}. The fluid motion is in the domain  $\Omega_a= [0,2\pi/a]\times [0,2\pi]$ and satisfies the spatially  periodic boundary condition \cite{CD}
 \begin{equation}\label{nb2}
\psi(2\pi/a,x_2)=\psi(0,x_2),\,\,\, \psi(x_1,0)=\psi(x_1,2\pi),\,\, \x=(x_1,x_2) \in {\Omega_a}.
\end{equation}
The  averaging condition
\begin{equation}\label{nbb3}
\int_{\Omega_a}  \psi dx_1dx_2 =0
\end{equation}
is applied to rule out non-zero constants being  solutions of the problem described by (\ref{nb1})--(\ref{nb2}).
  Note that $\psi=\psi^*$ is a steady-state
solution with respect to any $\kappa$. The solution multiplicity is thus obtained if there exists a family of
solutions $\psi_\kappa$ branching off $\psi^*$ from a critical value
$\kappa_{\rm{crit}}>0$.

The main result of the present paper reads as follows:
\begin{Theorem} \label{main}
 For  $1/\sqrt{2}\leq a<1$,
the equations  (\ref{nb1})--(\ref{nbb3}) admit a  positive critical value
\bbe\label{eq20}
\kappa_a< a\sqrt{\frac{1-a^2}{2(1+a^2)}}
\bee
  and a continuous family of classical steady-state solutions
$(\psi_{\kappa}, \kappa)$  branching off the bifurcation point $(\psi^*, \kappa_a)$ when $\kappa$ varies across  $\kappa_{a}$.
\end{Theorem}


 This result shows  mechanism behind the existence of a basic steady-state solution bifurcating into two  steady-state solutions under the single zonal  forcing (\ref{n3}). With the thermal forcing (\ref{n3}), the small $\kappa$ value implies that
 the acceleration nonlinearity dominates the circulation flow and then gives rise to multiple steady-state solutions, whereas the increment of the $\kappa$ value enlarges the linear Ekman layer dissipation and then eventually eliminates  the bifurcation phenomenon.

 Thus (\ref{nb1})  is quite similar to Navier-Stokes equations that the Ekman force $\kappa \Delta\psi$ plays  the same roll as the Reynolds viscous force $\frac1{\rm{Re}}\Delta^2 \psi$ to control the solution uniqueness and bifurcation behaviours.
 For the connection to the Euler equations, the Ekman dissipation force $\kappa \Delta \psi$ was recently unitized by the author \cite{C1,C2} to form a dissipative potential flow and then to produce dissipative free-surface Green functions for the cancelation of  wave integral singularity in numerical simulations of body motions in free water waves.

The equation (\ref{nb1}) is a third-order strong nonlinear  partial differential equation and is quite different to traditional semilinear  fluid motion equations such as the Navier--Stokes equations discussed in Temam \cite{Te} and the quasi-geostrophic equations discussed in Chen et al. \cite{Chen} and Chen and Price \cite{Chen1}. The  semilinearity indicates that the non-linear    term can be controlled by the  linear term. Therefore the a priori estimates and compactness analysis  of Navier--Stokes type equations,  available due to the presence of viscous force  (see, for example,  \cite{Chen,Chen1,Te}),  are not applicable to the strong non-linear   equation  (\ref{nb1}).
Actually, the non-linear  term of (\ref{nb1}) is the total derivative of fluid velocity along a particle trajectory and hence it is beneficial to use the Lagrangian formulation instead of the Euler formulation (\ref{nb1}) to control the nonlinearity of (\ref{nb1}).

For the equation (\ref{nb1}) with  the Dirichlet boundary condition, when the external forcing is changed into multiple ones  the existence of multiple steady-state responses    was discussed by Wolansky \cite{Wol2}.
The present state-state bifurcation analysis is applicable to the Dirichlet boundary value problem.  However, for the  vorticity  equation driven by a single forcing, it was unknown whether  the basic  solution branches into multiple steady-state  solutions when the Ekman dissipation force varies. Moreover  the steady-state bifurcation  analysis of the present paper, using the Krasnoselskii bifurcation theorem \cite{Kra} and the linear spectral technique developed from Meshalkin and Sinai \cite{MS}, Iudovich \cite{Iu} and Chen et al.  \cite{Chen} and Chen and Price \cite{Chen1}, is quite different to the multiple solution technique of Wolansky \cite{Wol2} although  the Lagrangian formulation  is developed from Wolansky \cite{Wol1}.

The functions in the present paper are in the H\"older spaces $C^{k+\alpha}({\Omega_a})$ for integer $k\geq 0$ and real
$\alpha \in [0,1)$. Here $C^0({\Omega_a})$ is the Banach space of all continuous functions  over ${\Omega_a}$ under the norm $$\|\phi\|_{C^0}=\max_{\x \in {\Omega_a}}|\phi(\x)|.$$
The $C^k$ and $C^{k+\alpha}$ function spaces are defined as $$C^k({\Omega_a})=\{ \phi \in C^0({\Omega_a}); \,\, \nabla^k\phi \in C^0({\Omega_a})\} $$
 with  the norm $$ \|\phi\|_{C^k}=\|\phi\|_{C^0}+\|\nabla^k \phi\|_{C^0},$$
$$C^{k+\alpha}({\Omega_a})=\{ \phi \in C^k({\Omega_a});\,\,\, \|\phi\|_{C^{k+\alpha}}=\|\phi\|_{C^k}+[\nabla^k \phi]_{C^\alpha}\},\,\,\, 0<\alpha<0,
$$
with the semi-norm
$$
[\phi]_{C^\alpha}=\sup_{\x,\y \in {\Omega_a}, \,\x\neq \y}\frac{|\phi(\x)-\phi(\y)|}{|\x-\y|^\alpha}.$$
 We use the function space $$C^{k+\alpha}_{\rm{per}}({\Omega_a})=\left\{ \phi \in C^{k+\alpha }({\Omega_a} );\,\, \mbox{ $\phi$ satisfies the conditions (\ref{nb2}) and (\ref{nbb3})}\right\}.
$$
A  steady-state solution $\psi$ of   (\ref{nb1})--(\ref{nbb3})  is said to be regular if
$\psi \in C^2_{\rm{per}}(\Omega_a)$ and $\Delta \psi \in C^1_{\rm{per}}(\Omega_a )$.

This paper is organized as follows. Section 2 exhibits a Lagrangian formulation approach to
 the atmospheric  flow in a neighborhood of the basic flow $\psi^*$ so that the compactness required by  the bifurcation analysis is obtained. Section 3 is devoted to the linear spectral analysis of the  vorticity   equation in the Lagrangian formulation. The spectral analysis technique is essentially developed from \cite{Chen,Chen1,Iu,MS}.
With the preparations of the compactness and the spectral results, Section 4 is devoted to the verification of the conditions ensuring the occurrence of the steady-state bifurcation phenomenon in  Krasnoselskii's  theorem. The proof of  Theorem \ref{main} is finally  completed in Section 4.

\section{Lagrangian formulation of the fluid motion}

For the velocity $\u =(u_1,u_2)=\nabla \times \psi$ of the fluid flow in the domain $\Omega_a$ and a  trajectory $\y=(y_1,y_2)$ initiating from a particle $\x=(x_1,x_2)$, the fluid motion is described by the Lagrangian formulation
\begin{equation}\label{qflow}\left\{ \begin{array}{rll}-\displaystyle\frac{\p}{\p t} \y(\x,t)&=
\u(\y(\x,t)),& t>0,
\vspace{2mm}
\\
\y(\x,0)&=\x \in {\Omega_a}. &
\end{array}\right.
\end{equation}
 Thus for the operators $$\nabla =(\partial_{x_1},\p_{x_2}),\,\,\,\,\,\nabla_\y= (\p_{y_1},\p_{y_2}),\,\,\,\,\nabla \y\cdot \nabla_\y =(\nabla y_1)\p_{y_1}+(\nabla y_2)\p_{y_2}$$ and the   $2\times 2$ identity  matrix $I$, we have
\bbe \label{yy}-\frac{\p }{\p t}\nabla\y(\x,t)
&=&\nabla \y\cdot\nabla_\y \u(\y(\x,t)),\,\,\, t>0,
\vspace{2mm}
\\
\nabla \y(\x,0)&=&I.\label{new5}
\bee
This system  implies the Euler identity
\be -\frac{\p}{\p t} \det (\nabla \y)= \det (\nabla \y)\nabla_\y \cdot \u(\y)
\ee
and hence the incompressible flow transformation property
\bbe \label{sss}
\det(\nabla \y)=1.
  \bee
It follows from (\ref{yy}) that
\bbe\label{new4}
\frac{\p }{\p t} |\p_{x_i} y_j| \leq |\p_{x_i} y_1|\,|\p_{y_1}u_j(\y)|+|\p_{x_i} y_2|\,|\p_{y_2}u_j(\y)|,\,\,\, i,\,j=1,\,2.
\bee
Here the time derivative  $\p_t |f|$ is in the sense of $ \limsup_{\delta t \to 0}\frac{|f(t+\delta t)|-|f(t)|}{\delta t}$.
We thus have
\be \lefteqn{\frac12 \frac{\p }{\p t}|\nabla \y|^2}\\
&\leq& \left(|\p_{x_1} y_1|^2+|\p_{x_2} y_1|^2\right)|\p_{y_1}u_1(\y)|+\left(|\p_{x_1} y_2|^2+|\p_{x_2} y_2|^2\right)\,|\p_{y_2}u_2(\y)|
\\
&&+\frac12\left(|\p_{x_1} y_2|^2+|\p_{x_1} y_1|^2+|\p_{x_2} y_2|^2+|\p_{x_2} y_1|^2\right)\,\left(|\p_{y_1}u_2(\y)|+|\p_{y_2}u_1(\y)|\right)
\\
&\leq& \frac 54 |\nabla \y|^2 \|\nabla_\y \u\|_{C^0}.
\ee
This together with (\ref{new5}) gives the flow estimate  expressed as
\bbe \label{new3}
 |\nabla \y(\x,t)| &\leq& \sqrt{2}\e^{\frac 54 t\|\nabla_\y \u\|_{C^0}}.
\bee

On the other hand,
 the study of the uniqueness and the multiplicity of the classical solutions around the basic solution $\psi^*$
    is based on the flow estimate  expressed as
 \bbe\label{aaaa2} |\nabla \y(\x,t)| \leq  (\sqrt{2}+\sqrt{5}t)\e^{2t  \sqrt{\|\nabla_\y\u-\nabla_\y \u^*\|_{C^0}}}\bee
for $\u^*=\nabla \times \psi^*$. Hence, for convenience, we may assume that the inequality
\bbe\label{new2}\|\nabla_\y\u-\nabla_\y \u^*\|_{C^0}\leq \frac12
\bee is always  true  since
the present investigation aims at the  uniqueness and bifurcation around the basic flow $\psi^*$.

To show the validity of (\ref{aaaa2}), we set $\epsilon =\|\nabla_\y\u-\nabla _\y\u^*\|_{C^0}$ or
\be \epsilon = \left\|\sqrt{(\p_{y_1}u_1)^2+(\p_{y_2}u_1-\cos y_2)^2+(\p_{y_1}u_2)^2+(\p_{y_2}u_2)^2}\right\|_{C^0}.\ee
With the use of
the matrix inequality notation
$$(a_{i,j}) \leq (b_{i,j}) \mbox{ whenever } a_{i,j}\leq b_{i,j} \mbox{ for all $i$ and $j$,}
$$
 the equation (\ref{new4}) can be rewritten as
\be
\frac{\p }{\p t}\left(\begin{array}{ll}|\p_{x_1}y_1|&|\p_{x_2}y_1|\vspace{1mm}\\
|\p_{x_1}y_2|&|\p_{x_2}y_2|
\end{array}\right)&\leq& \left(\begin{array}{ll}|\p_{y_1}u_1|&|\p_{y_2}u_1|\vspace{1mm}\\
|\p_{y_1}u_2|&|\p_{y_2}u_2|
\end{array}\right)\left(\begin{array}{ll}|\p_{x_1}y_1|&|\p_{x_2}y_1|\vspace{1mm}\\
|\p_{x_1}y_2|&|\p_{x_2}y_2|
\end{array}\right)
\\
&\leq & \left(\begin{array}{lll}\epsilon &&1+\epsilon \vspace{1mm}\\
\epsilon &&\epsilon
\end{array}\right)\left(\begin{array}{ll}|\p_{x_1}y_1|&|\p_{x_2}y_1|\vspace{1mm}\\
|\p_{x_1}y_2|&|\p_{x_2}y_2|
\end{array}\right)
\ee

Multiplying this inequality by the matrix
$$ \exp \left( -t\left(\begin{array}{lll}\epsilon &&1+\epsilon \vspace{1mm}\\
\epsilon &&\epsilon
\end{array}\right)\right)
$$
and using  the initial condition $\nabla \y(\x,0) =I$, we have
\be \lefteqn{
\left(\begin{array}{ll}|\p_{x_1}y_1|&|\p_{x_2}y_1|\vspace{1mm}\\
|\p_{x_1}y_2|&|\p_{x_2}y_2|
\end{array}\right)}
\\
&\leq& \exp \left( t\left(\begin{array}{lll}\epsilon &&1+\epsilon \vspace{1mm}\\
\epsilon &&\epsilon
\end{array}\right)\right)
\\
&=&\left(\begin{array}{lr}
\frac {\sqrt{\epsilon^2+\epsilon}}{2\epsilon}&\frac {\sqrt{\epsilon^2+\epsilon}}{2\epsilon}\vspace{2mm}
\\
-\frac12 & \frac12
\end{array}\right)\left(\begin{array}{lr}\e^{t(\epsilon +\sqrt{\epsilon^2+\epsilon})} &0\vspace{2mm}\\
0&\e^{t(\epsilon -\sqrt{\epsilon^2+\epsilon})}
\end{array}\right)\left(\begin{array}{lr}\frac\epsilon {\sqrt{\epsilon^2+\epsilon}} &\,\,\,\,-1\vspace{2mm}\\
\frac\epsilon {\sqrt{\epsilon^2+\epsilon}}&1
\end{array}\right)
\ee
and hence, for $\lambda_1=\epsilon +\sqrt{\epsilon^2+\epsilon}$ and  $\lambda_2=\epsilon -\sqrt{\epsilon^2+\epsilon}$,
\be |\nabla \y(\x,t)|^2&\leq&   \frac{\left(\e^{\lambda_1t}+\e^{\lambda_2 t}\right)^2
+ \left(\frac{\epsilon +1}\epsilon +\frac\epsilon{\epsilon +1}\right)
\left(\e^{\lambda_1 t}-\e^{\lambda_2 t}\right)^2}2
\\
&\leq&   \left[2
+ \frac12\left(\frac{\epsilon +1}\epsilon +\frac\epsilon{\epsilon +1}\right)
\left(\lambda_1t-\lambda_2t\right)^2\right]\e^{2t(\epsilon +\sqrt{\epsilon^2+\epsilon})}
\\
&\leq&   \left[2
+ 2t^2[(\epsilon +1)^2+\epsilon^2] \right]\e^{2t(\epsilon +\sqrt{\epsilon^2+\epsilon})}
\\
&\leq& (2+5t^2)\e^{2t(\epsilon +\sqrt{\epsilon^2+\epsilon}) }.
\ee
Here we have  used  equation (\ref{new2}).
The validity of    (\ref{aaaa2}) is thus demonstrated.

The following lemma shows the well-posedness of the fluid motion in the Lagrangian formulation:
\begin{Lemma}\label{alpha}
 Assume that    $\kappa>0$ and   $\psi \in C^{2}_{\rm{per}}({\Omega_a})$ such that
   \bbe \label{yyy} \label{asas}
  \|\nabla^2\psi-\nabla^2\psi^*\|_{C^0}< \frac{\kappa^2}4.
  \bee
  Then
 the operator
$\kappa+(\nabla\times \psi)\cdot \nabla $ is a bijection mapping the space
$$D= \{ f \in
C^{1}_{\rm{per}}({\Omega_a});\,\,(\kappa+(\nabla\times \psi)\cdot \nabla)   f \in
C^{1}_{\rm{per}}({\Omega_a})\}$$
onto $C^{1}_{\rm{per}}({\Omega_a})$ and
\be
\lefteqn{\|[\kappa+(\nabla\times \psi)\cdot \nabla]^{-1}f\|_{C^{1}}}
\\
&\leq&\left(
\displaystyle\frac{1}\kappa+\displaystyle\frac{\sqrt{2}(\kappa-2\sqrt{\|\nabla^2\psi-\nabla^2\psi^*\|_{C^0}})+\sqrt{5}}
{(\kappa-2\sqrt{\|\nabla^2\psi-\nabla^2\psi^*\|_{C^0}})^2} \right)\|f\|_{C^{1}}.
\ee
\end{Lemma}

\noindent {\sl Proof}.
For the injection assertion, we see that the equation
$$(\kappa+(\nabla\times \psi)\cdot \nabla )  f =0$$
 implies, with the use of integration by parts,
\begin{eqnarray*}
\kappa\int_{{\Omega_a}} f ^2dx_1dx_2&=&-\int_{{\Omega_a}} f (\nabla\times \psi)\cdot \nabla  fdx_1dx_2 \\
&=&\int_{{\Omega_a}}  f (\nabla\times \psi) \cdot \nabla fdx_1dx_2 \\
&=&-\kappa\int_{{\Omega_a}} f ^2 dx_1dx_2,
\end{eqnarray*}
which shows  $f=0$.

For the surjection assertion,   we consult   \cite{Wol1} to  define the operator
 $$T_\psi f(\x)=\int^\infty_0
e^{-\kappa s} f (\y(\x,s))ds,$$
which is utilized  to show the required conditions
$$T_\psi f\in D\,\,\,\mbox{ and }\,\,\,
(\kappa+(\nabla\times \psi)\cdot \nabla )  T_\psi f= f.$$

Indeed, upon the observation of the equation
\bbe\label{psif}\nabla T_\psi f(\x) =\int^\infty_0
e^{-\kappa s}\nabla \y(\x,s)\cdot \nabla_\y  f (\y(\x,s)) ds
 \bee
 and the quantity  $\epsilon = \|\nabla ^2 \psi-\nabla^2\psi^*\|_{C^0}$, it follows from  (\ref{aaaa2}) that
 \begin{eqnarray*}|T_\psi f (\x)|+|\nabla T_\psi f(\x)
  |&\leq &\int^\infty_0
e^{-\kappa s}\|  f \|_{C^0}ds + \int^\infty_0 e^{-\kappa s}\|\nabla
 f \|_{C^0}\|\nabla \y\|_{C^0}ds\\ &\leq
&\frac1\kappa\| f \|_{C^0}+\int^\infty_0(\sqrt{2}+\sqrt{5}s) e^{ - \kappa s+2\sqrt{\epsilon} s}
\|\nabla  f \|_{C^0} ds\\ &\leq
&\left(\frac1\kappa+\frac{\sqrt{2}}{\kappa-2\sqrt{\epsilon}}+
\frac{\sqrt{5}}{(\kappa-2\sqrt{\epsilon})^2}\right)\| f \|_{C^1},
\end{eqnarray*}
which gives the estimate of the operator $T_\psi$.

To verify the continuity of the function $\nabla T_\psi f$,
we employ (\ref{yy}) and (\ref{aaaa2}) to produce
\bbe \label{y}| \y(\x,t)-\y(\x',t)| \leq (\sqrt{2}+\sqrt{5}t) e^{ 2\sqrt{\epsilon} t}|\x-\x'|,\,\,\,\,\x,\,\x'\in {\Omega_a},
\bee
and
\be -\frac{\p }{\p t}(\nabla\y(\x,t)-\nabla \y(\x',t))
&=&(\nabla \y(\x,t)-\nabla \y(\x',t))\cdot\nabla_\y \u(\y(\x,t))
\\
&&+ \nabla \y(\x',t)\cdot (\nabla_\y \u(\y(\x,t))-\nabla_\y \u(\y(\x',t))).
 \ee
 Hence the derivation of (\ref{aaaa2}) implies
  \bbe \lefteqn{|\nabla\y(\x,t)-\nabla \y(\x',t))|}\nonumber
\\
\nonumber
&\leq &\int^t_0 (\sqrt{2}+\sqrt{5}t) e^{ 2\sqrt{\epsilon} (t-s)}| \nabla \y(\x',s)|\,|\nabla_\y \u(\y(\x,s))-\nabla_\y \u(\y(\x',s))|ds
\\ &\leq &\e^{2t\sqrt{\epsilon}}(\sqrt{2}+\sqrt{5}t)^2\int^t_0  |\nabla_\y \u(\y(\x,s))-\nabla_\y \u(\y(\x',s))|ds.\label{dy}
 \bee

Moreover, for any  constant $\tau>0$, it follows  from (\ref{aaaa2}) and (\ref{psif}) that
\bbe \nonumber
\lefteqn{|\nabla T_\psi f (\x)-\nabla T_\psi f(\x')|}\\
  &\leq &\int^\infty_0\nonumber
e^{-\kappa s+2s\sqrt{\epsilon}}(\sqrt{2}+\sqrt{5}s)|\nabla_\y f (\y(\x,s))-\nabla _\y f (\y(\x',s))|ds\\
&& +\| \nabla f\|_{C^0} \int^\infty_0 e^{-\kappa s}\, |\nabla
\y(\x,s)-\nabla \y(\x',s)|ds\nonumber
\\
&\leq &3\|\nabla f\|_{C^0}\int^\infty_\tau
e^{-\kappa s+2s\sqrt{\epsilon}}(\sqrt{2}+\sqrt{5}s)ds\label{df1}
\\
&&+\int^\tau_0
e^{-\kappa s+2s\sqrt{\epsilon}}(\sqrt{2}+\sqrt{5}s)|\nabla_\y f (\y(\x,s))-\nabla _\y f (\y(\x',s))|ds\label{df2}
\\
&&+\| \nabla f\|_{C^0} \int^\tau_0 e^{-\kappa s}\, |\nabla
\y(\x,s)-\nabla \y(\x',s)|ds.\label{df3}
\bee
Therefore, for any $\varepsilon>0$, we can use (\ref{aaaa2}), (\ref{yyy}), (\ref{y}) , (\ref{dy}) and the continuity of $\nabla f$ and $\nabla \u$   to demonstrate that each of the items (\ref{df1})--(\ref{df3}) is bounded by $ \varepsilon/3$, provided that $\tau>0$ is sufficiently large  and
$|\x-\x'|$ is sufficiently small.
Hence $T_\psi f\in C^1_{\rm{per}}({\Omega_a})$.

  The surjection is due to the validity of the identity
 \be (\kappa+(\nabla \times \psi)\cdot \nabla )T_\psi f= f,\ee
which is demonstrated as follows:
\begin{eqnarray*}
(\nabla\times \psi)\cdot \nabla
T_\psi f(\x)
 &=&\lim_{t\to 0+}\left.
(\nabla_\y\times \psi(\y(\x,t)))\cdot \nabla_\y
T_\psi f(\y(\x, t ))\,\right.
\\ &=&-\lim_{t\to 0+}\left.
\frac{\partial \y(\x, t )}{\p t}  \cdot
\nabla_\y T_\psi f(\y(\x, t ))\,\right.,\,\,\mbox{ by (\ref{qflow})},
\\
&=&-\lim_{t\to 0+}\left. \frac{\p}{\p t}
T_\psi f(\y(\x, t ))\,\right.
\\ &=&\left. -\lim_{t\to 0+}\frac{\p}{\p t} \int^\infty_0e^{- \kappa s} f (\y(\x, t
+ s))ds \right.
\\ &=&-\lim_{t\to 0+}\left.
\frac{\p}{\p t} \int_ t ^{\infty}e^{  -\kappa (s-t)} f (\y(\x, s))ds
\right.=-\kappa T_\psi f (\x) + f(\x) .
\end{eqnarray*}
The proof is completed.

As a consequence of Lemma \ref{alpha}, the steady-state problem of the Euler formulation (\ref{nb1})--(\ref{nbb3}) becomes the
Lagrangian formulation problem
\be-\Delta  \psi = \kappa [\kappa +(\nabla\times\psi)\cdot \nabla ]^{-1}\psi^*\ee
or
\bbe  \label{Lag}-\Delta \psi(\x)=\kappa\int^\infty_0 \e^{-\kappa s} \psi^*(\y(\x,s))ds,
\bee
provided that $\psi\in C^2_{\rm{per}}({\Omega_a})$ satisfies  the condition (\ref{asas}).

It is readily seen that the proof of Lemma \ref{alpha} remains true if we utilize the estimate (\ref{new3})  instead of  the estimate (\ref{aaaa2}). More precisely,
  the proof of Lemma \ref{alpha}  implies    the following  regularity criterion.
\begin{Lemma}\label{cor}
 For $0<\alpha <1$ and $\kappa>0$, let   $\psi\in C^2_{\rm{per}}({\Omega_a})$ be  a solution of  (\ref{Lag}) satisfying either the condition (\ref{asas}) or the condition
 \bbe\label{new1} \|\nabla^2 \psi\|_{C^0}  <\frac 45\kappa.\bee
 Then $\psi \in C^{2+\alpha}_{\rm{per}}({\Omega_a})$ and  $\Delta \psi \in C^1_{\rm{per}}(\Omega_a)$. That is, $\psi$  is  a regular  solution of the problem described by (\ref{nb1})--(\ref{nbb3}).
 \end{Lemma}

The uniqueness assertion of  Theorem \ref{main}  is implied from the following.
\begin{Theorem}\label{p2}

Let  $\kappa\geq \frac1a$ and $\psi\in C^2_{\rm{per}}(\Omega_a)$ be  a solution of  the Lagrange formulation problem (\ref{Lag}) or the Euler formulation problem (\ref{nb1})--(\ref{nbb3}) satisfying
the condition (\ref{asas}). Then $\psi$ is regular  and
$\psi=\psi^*$ holds true.\label{cor2}
\end{Theorem}

The uniqueness  was  discussed  \cite{Titi1988,Hauk1997,Il1993,Wol1} in the vicinity of a small steady-state solution. In contrast,  Theorem \ref{p2} is on  the uniqueness   in the vicinity of the basic steady-state  solution $\psi^*$, which is not small.

\

\noindent
{\sl Proof}.  We employ  Lemma  \ref{cor}  to obtain the regularity of  $\psi$, which is a steady-state  solution of (\ref{nb1})--(\ref{nbb3}). The observation
\bbe \label{psip}-\Delta \psi^*=\kappa [\kappa+(\nabla\times \psi^*)\cdot\nabla ]^{-1}\psi^*
\bee
and  the  application of  the $L_2$ norm
$$\|\phi\|_{L_2}=\left(\int_{\Omega_a} |\phi(\x)|^2dx_1dx_2\right)^{1/2}$$
yield that
\be
\lefteqn{\|\Delta \psi-\Delta \psi^*\|_{L_2}}\\&=&\kappa\|[\kappa+(\nabla\times \psi)\cdot\nabla ]^{-1}\psi^*-[\kappa+(\nabla \times \psi^*)\cdot\nabla ]^{-1}\psi^*\|_{L_2}
\\&=&\kappa \|[\kappa+(\nabla\times \psi)\cdot\nabla ]^{-1}(\nabla \times \psi-\nabla\times \psi^*)\cdot \nabla [\kappa+(\nabla\times \psi^*)\cdot\nabla ]^{-1}\psi^*\|_{L_2}
\\&\leq& \|(\nabla \times \psi-\nabla\times \psi^*)\cdot \nabla [\kappa+(\nabla\times \psi^*)\cdot\nabla ]^{-1}\psi^*\|_{L_2},
\ee
where we have used the variable transformation property (\ref{sss}) and the integral formulation  (\ref{Lag}). By  (\ref{psip}), we thus have
\be
\|\Delta \psi-\Delta \psi^*\|_{L_2}&\leq &\frac1\kappa \|(\nabla\times \psi-\nabla\times \psi^*)\cdot \nabla \psi^*\|_{L_2}
\\&\le&\frac1\kappa \|\nabla\times \Delta ^{-1}\Delta (\psi- \psi^*)\|_{L_2}
\\&< &\frac1{\kappa a} \|\Delta \psi- \Delta \psi^*\|_{L_2}.
\ee
whenever $\psi\neq \psi^*$.  This leads to a contradiction since $\kappa a \geq 1$. Hence  $\psi=\psi^*$.

The proof of Theorem \ref{p2} and hence the proof  of Theorem \ref{main} (i) are   completed.

\section{Linear spectral analysis}
For   steady-state  solutions  branching off the basic  solution $\psi^*=\cos x_2$ or the existence of  steady-state solutions in a vicinity  of $\psi^*$, it follows from Lemma \ref{alpha} that the steady-state  Euler formulation  problem
(\ref{nb1})--(\ref{nbb3}) is equivalent to the Lagrangian formulation problem
\bbe \label{bh1}\psi +\kappa \Delta ^{-1}[\kappa +(\nabla \times \psi)\cdot \nabla ]^{-1} \psi^*=0,\,\,\,\, \psi\in C^2_{\rm{per}}({\Omega_a}),\,\, \Delta \psi \in C^1(\Omega_a).\bee
However the bifurcation phenomenon of  (\ref{bh1}) results from the nonlinearity and linear spectral analysis  of the problem (\ref{bh1}). This section is contributed to the spectral analysis of the
 operator $L_\kappa$ linearized from  the non-linear  operator $F(\psi,\kappa)$, the left-hand side term of  (\ref{bh1}), around the basic flow  $\psi^*$. By an elementary manipulation, the operator $L_\kappa$ can be linearized  as
 \bbe\nonumber
L_\kappa\psi&=&\displaystyle\lim_{s\to 0} \frac{F(\psi^*+s\psi,\kappa)-F(\psi^*,\kappa)}{s}
 \\ \nonumber &=&\psi +\lim_{s\to 0}\left( \frac{ \kappa\Delta^{-1}[\kappa +(\nabla \times \psi^*+s\nabla \times \psi)\cdot \nabla ]^{-1}\psi^*}s
 \right.
 \\
 &&-\left.\frac{\kappa\Delta^{-1}[\kappa+(\nabla\times \psi^*)\cdot \nabla ]^{-1}\psi^*}{s}\right)\nonumber
 \\
 &=&\psi - \kappa\Delta^{-1}[\kappa +(\nabla \times \psi^*)\cdot \nabla ]^{-1}(\nabla \times \psi)\cdot \nabla [\kappa+(\nabla\times \psi^*)\cdot \nabla ]^{-1}\psi^*\nonumber
 \\&=&\psi - \Delta^{-1}[\kappa +(\nabla \times \psi^*)\cdot \nabla ]^{-1}(\nabla \times \psi)\cdot \nabla \psi^*\nonumber
   \\&=&\psi + \Delta^{-1}[\kappa +\sin x_2 \p_{x_1} ]^{-1}(\sin x_2\p_{x_1} \psi),\label{Lk}
 \bee
 where we have used the solution property (\ref{psip}).

We can now examine the critical  real spectral problem
\bbe \label{eigen}L_\kappa \psi = 0
\bee
in the space $C^{2+\alpha}_{\rm{per}}({\Omega_a})$ with $0\leq \alpha \leq 1$.
Here $\kappa$ is said to be a critical  if
equation (\ref{eigen})  admits a non-zero solution or  an eigenfunction $\psi\in
C^{2+\alpha}_{\rm{per}}({\Omega_a})$.
The spectral problem is restricted in the even function subspace
\begin{eqnarray*}\hat
C^{2+\alpha}_{\rm{per}}({\Omega_a})&=&\left\{\psi\in  C^{2+\alpha}_{\rm{per}}({\Omega_a});\,\,
\psi(-\x)=\psi(\x)\right\}.
\end{eqnarray*}
By Fourier expansion, the function $\psi$ in $\hat
C^{2+\alpha}_{\rm{per}}({\Omega_a})$ is generally expressed as
$$
\psi=\sum_{m=0}^\infty \sum_{n=-\infty}^\infty b_{m,n}\cos(max_1+nx_2).
$$

The spectral result is stated as follows:
\begin{Theorem} \label{spectral} Let $\frac1{\sqrt{2}}\leq a<1$ and $\kappa >0$.
Then  there exists a positive critical value
 $$\kappa_a< a\sqrt{\frac{1-a^2}{2(1+a^2)}}$$ such that
\bbe\label{one11}\dim\bigcup_{i=1}^\infty\left\{\psi=\sum_{n=-\infty}^\infty b_n\cos (a x_1+nx_2) \in \hat C^{2+\alpha}_{\rm{per}}({\Omega_a});\,\,L_{\kappa_a}^i \psi=0\right\}= 1.\bee
If $m\neq 1$ is a nonnegative  integer, then it is valid that
\bbe \label{zero}\dim\left\{\psi=\sum_{n=-\infty}^\infty b_n\cos (ma x_1+nx_2) \in \hat C^{2+\alpha}_{\rm{per}}({\Omega_a});\,\,L_\kappa \psi=0\right\}= 0.
\bee
\end{Theorem}
Theorem \ref{spectral} is   proved  by a  continued fraction technique  developed from Chen et al.  \cite{Chen} and  Chen and Price \cite{Chen1} and originated from
  Mishalkin and Sinai \cite{MS} and  Iudovich
\cite{Iu}.
However,  the linear operator $L_\kappa$ now involves  the Lagrangian formulation aspect.

\

\noindent {\sl Proof}. To verify the validity of (\ref{zero}), we use (\ref{Lk}) to rewrite the spectral equation $L_\kappa\psi=0$ as
 \bbe \label{new0}\Delta \psi + (\kappa +\sin x_2 \p_{x_1} )^{-1}(\sin x_2\p_{x_1} \psi)=0.
 \bee
 It is readily seen that the operator
 $$(\kappa +\sin x_2 \p_{x_1} )^{-1}=[\kappa+(\nabla \times \psi^*)\cdot \nabla ]^{-1}$$ maps $C^1$ into $C^1$ or $\Delta \psi \in C^1$.
 Thus we may apply  the operator $(\kappa+\sin
 x_2\p_{x_1})$ to (\ref{new0}) to produce  the spectral equation
\begin{equation}\kappa\Delta\psi + \sin x_2 (\Delta +1)\p_{x_1} \psi =0.\label{eigen1} \end{equation}
Multiplying (\ref{eigen1}) by $(\Delta +1)\psi$ and integrating the resultant equation over the domain ${\Omega_a}$, we have
the integral equation
\bbe\label{stable} 0=\int_{{\Omega_a}} \Delta \psi (\Delta \psi+\psi) dx_1dx_2.
\bee
The substitution of the function $$ \psi =\sum_{n=-\infty}^\infty
b_{n}\cos(amx_1+nx_2),\,\,\, m\neq 1.
$$
into  (\ref{stable}) simply implies  $b_n\equiv 0$ and hence  (\ref{zero}) is verified. Here we have used the  average condition (\ref{nbb3}) to confirm  $b_0=0$ whenever $m=0$.

To show the existence of the critical number $\kappa_a$, we substitute the
eigenfunction
\begin{equation}\label{mn}
\psi =\sum_{n=-\infty}^\infty
b_{n}\cos(ax_1+nx_2)
\end{equation}
  into (\ref{eigen1})
to obtain the
 iteration
equation, for arbitrary integer $n$,
\bbe\label{diff0}
2\kappa (a^2+n^2)b_n-a[a^2+(n+1)^2-1]b_{n+1}+a[a^2+(n-1)^2-1]b_{n-1}=0
\bee
or
\bbe\label{diff}
\kappa d_n(\beta_n-1)b_n-(\beta_{n+1}-1)b_{n+1}+(\beta_{n-1}-1)b_{n-1}=0
\bee
for
\bbe\beta_n=a^2+n^2\,\mbox{ and }\, d_n
=\frac{2\beta_n}{a(\beta_n-1)}.\label{eq40}
\bee
 Notice  that $(\beta_n-1) b_n\neq 0$ for any
$n$ since $b_n \equiv 0$ if and only if $b_{n_0}=0$ for an integer $n_0$ (see \cite{MS}). This enables us to define the quantities
\bbe \label{eq30}\gamma_{n}= \frac{(\beta_n-1) b_n}{(\beta_{n-1}-1) b_{n-1}}, \,\,\, \gamma_{-n}=\frac{ (\beta_n-1) b_{-n}}{(\beta_{n-1}-1) b_{-n+1}} \,\,\,\mbox{ for }\,\,n > 0.
\bee
Thus by    dividing  (\ref{diff}) with the quntity $(\beta_n-1) b_n$,
the equation  (\ref{eq40}) is written as
\bbe kd_n - \gamma_{n+1} + \frac{1}{\gamma_n}=0\,\,\mbox{ for } n> 0,\label{new6}
\\
kd_n -  \frac{1}{\gamma_{-n}}+\gamma_{-n-1} =0\,\,\mbox{ for } n> 0,\label{new7}
\\
\label{eq8}d_0k-\gamma_1+\gamma_{-1} =0\,\,\mbox{ for } n=0.
\bee

 With the use of  (\ref{new6})--(\ref{new7}), we have
\begin{equation}\label{gamma}
\gamma_{\pm n} =\frac{\mp 1}{\kappa d_n\mp
\gamma_{\pm (n+1)} }=\frac{\mp
1}{\kappa d_n+\displaystyle\frac{1}{
\kappa d_{n+1}+\displaystyle\frac{1}{\ddots}}}\,\,\mbox{ for }\,\,n\geq
1.
\end{equation}
It follows from (\ref{eq40}), (\ref{eq8}) and (\ref{gamma}) that the spectral problem (\ref{eigen}) or (\ref{diff}) is
equivalent to the  equation
\begin{equation}\label{frac}
\frac{a}{1-a^2}=\frac{
1}{{\kappa^2}d_1+\displaystyle\frac{1}{
d_2+\displaystyle\frac{1}{
{\kappa^2}d_3+\displaystyle\frac{1}{
d_4+\displaystyle\frac{1}{ \ddots}}}}}.
\end{equation}

The function  $P(\kappa)$, representing  the right-hand side term  of  (\ref{frac}), is  the Stieltjes continued fraction. It follows from \cite{Sti} or \cite[Theorem
28.1]{Wal} that  $P(\kappa)$ uniformly convergent to a positive
value and is an analytic function of $\kappa >0$. Upon observation of  $P$  being
strictly monotone function of $\kappa$ such that
$$\lim_{\kappa\to
\infty}P(\kappa)=0,\,\,\,\,\lim_{\kappa\to 0}P(\kappa)=\infty,$$
there exists a unique critical value $\kappa=\kappa_a>0$
satisfying  (\ref{frac}). Thus for such a
critical value $\kappa=\kappa_a$, the coefficients $b_n$ of  the associated  eigenfunction $\psi$ in the form of  (\ref{mn}) and (\ref{eq30}) are subject to the
expression
\begin{equation}\label{oned}
b_n=\left\{\begin{array}{ll}\displaystyle c\frac{a^2-1}{a^2+n^2-1}\gamma_1\cdots
\gamma_n,& n\geq 1,\vspace{2mm}\\ c,&n=0,\vspace{2mm}\\
(-1)^n b_{-n},& n\leq -1.
\end{array}\right.\end{equation}
for an arbitrary  constant $c$. Equation (\ref{gamma}) implies that
\be \lim_{n\to \infty} \gamma_n = \frac{-1}{\displaystyle \frac{2\kappa}a -\displaystyle\lim_{n\to \infty} \gamma_n}
 \ee
 or
 \be  \lim_{n\to \infty} \gamma_n= \frac\kappa a-\sqrt{\frac{\kappa^2}{a^2}+1}=\frac{-1}{\frac ka +\sqrt{\frac{\kappa^2}{a^2}+1}}.\ee
 This gives the smoothness of the eigenfunction $\psi$  expressed by (\ref{mn}) and (\ref{oned}) and hence $\psi\in C^{2+\alpha}_{\rm{per}}({\Omega_a})$. That is,
 \bbe\label{one0}\dim\left\{\psi =\sum_{n=-\infty}^\infty b_n \cos(ax_1+n x_2)\in \hat C^{2+\alpha}_{\rm{per}}({\Omega_a});\,\,L_{\kappa_a}\psi=0 \right\}= 1.
\bee

The upper bound of the critical value $\kappa_a$ is an immediate consequence of the inequality
\be \frac{a}{1-a^2} \le \frac1{\kappa_a^2 d_1}=\frac{a^3}{2\kappa_a^2(a^2+1)},
\ee
which follows from  (\ref{eq40}) and (\ref{frac}).

To prove  the spectral simplicity given in  (\ref{one11}), it  is  sufficient to verify the property
\bbe\label{one}\dim\bigcup_{i=1}^2\left\{\psi =\sum_{n=-\infty}^\infty b_n \cos(ax_1+n x_2)\in \hat C^{2+\alpha}_{\rm{per}}({\Omega_a});\,\,L_{\kappa_a}^i\psi=0 \right\}= 1.
\bee

We see that the
 equation
$L_{\kappa_a}^2\psi=0$
can be written in the form
  \begin{equation}\label{con}L_{\kappa_a}\psi'=0\,\,\, \mbox{ and }\,\,\,\psi'=L_{\kappa_a}\psi\end{equation}
or, equivalently,
  \begin{eqnarray}\label{nnew1}
&&\kappa_a\Delta\psi' +\sin x_2 (\Delta +1)\p_{x_1} \psi' =0,
 \\
&&\kappa_a\Delta\psi + \sin x_2 (\Delta +1)\p_{x_1}\label{nnew2}
\psi=(\kappa_a + \sin x_2\p_{x_1})\Delta \psi'.
\end{eqnarray}
By the Fourier expansions
  $$\psi=\sum_{n=-\infty}^\infty
b_n\cos(ax_1+nx_2)\,\,\mbox{ and }\,\,\psi'=\sum_{n=-\infty}^\infty b_n'\cos(ax_1+nx_2),$$
 the equations (\ref{nnew1})-(\ref{nnew2}) reduce respectively  to the iteration equations
\bbe\label{eq1} 2\kappa_a\beta_n b'_n -
a(\beta_{n+1}-1)b'_{n+1}+a(\beta_{n-1}-1)b'_{n-1}=0
\bee
and
\bbe
&&2\kappa_a \beta_n b_n -a(\beta_{n+1}-1)b_{n+1}+
a(\beta_{n-1}-1)b_{n-1}\nonumber
\\ \label{eq2}
&&=2\kappa_a \beta_nb_n'-a \beta_{n+1} b_{n+1}'+a\beta_{n-1}b_{n-1}'
\bee
 for any arbitrary integer $n$.
Therefore from the demonstration of the assertion (\ref{one0}), it remains to prove  that $\psi'=0$ or $b_n'\equiv 0$.  Due to the   equivalence of   (\ref{diff}) and  (\ref{eq1}), all the equations involving the proof of (\ref{one0}) hold true if $b_n$ is replaced by $b_n'$ therein.

 Multiplying the $n$th equation of (\ref{eq1}) by $(-1)^n(\beta_n-1)b_n$ and the $n$th equation of (\ref{eq2}) by $(-1)^n(\beta_n-1)b_n'$ and then summing the resultant equations respectively, we have
\begin{eqnarray}\label{eq3}
0=\sum_{n=-\infty}^\infty (-1)^n(\beta_n\!-\!1)b_n[2\kappa_a\beta_n b_n'
-a (\beta_{n+1}\!-\!1)b_{n+1}'\!+\!a
(\beta_{n-1}\!-\!1)b_{n-1}']
\bee
and
\bbe\nonumber\lefteqn{ \sum_{n=-\infty}^\infty
(-1)^n(\beta_n-1)b_n'[2\kappa_a\beta_n b_n -a
(\beta_{n+1}-1)b_{n+1}+a (\beta_{n-1}-1)b_{n-1}]}
\\ &=&\label{eq4}
\sum_{n=-\infty}^\infty (-1)^n(\beta_n-1)b_n'[2\kappa_a\beta_n b_n' -a\beta_{n+1}b_{n+1}'+a\beta_{n-1}b_{n-1}')].
\end{eqnarray}
Rearranging terms in the summations, we see that the right-hand side term of  (\ref{eq3}) is identical to the  left-hand side term of (\ref{eq4}). Thus   (\ref{eq4}) becomes
\bbe \label{eq5}
\hspace{-3mm}0= \!\!\!\sum_{n=-\infty}^\infty \! \!\!(-1)^n2\kappa_a(\beta_n\!-\!1)\beta_nb_n'^2\!-\!\!\!\!\sum_{n=-\infty}^\infty\!\! \! (-1)^na(\beta_n\!-\!1)b_n'[\beta_{n+1}b_{n+1}'\!-\!\beta_{n-1}b_{n-1}'].
\end{eqnarray}
Therefore  it remains to show that (\ref{eq5}) leads to $b_n'\equiv 0$.
To do so, we formulate  the second term on the right-hand side of  (\ref{eq5}) as follows:
\be\nonumber
\lefteqn{ \sum_{n=-\infty}^\infty (-1)^na(\beta_n-1)b_n'[
\beta_{n+1}b'_{n+1}- \beta_{n-1}b'_{n-1}]}
 \nonumber
\\&=&\sum_{n=0}^\infty a(-1)^n(\beta_n-1)\beta_{n+1}b_n'b_{n+1}'
+\sum_{n=1}^{\infty }a(-1)^n(\beta_{-n}-1)\beta_{-n+1}b_{-n}'b_{-n+1}'
\\ \nonumber
&&-\sum_{n=1}^\infty a(-1)^n(\beta_n-1)\beta_{n-1}b_n'b_{n-1}'
-\sum_{n=0}^{\infty }a(-1)^n(\beta_{-n}-1)\beta_{-n-1}b_{-n}'b_{-n-1}'
\\ 
&=&\sum_{n=0}^\infty 2a(-1)^n(\beta_n-1)\beta_{n+1}b_n'b_{n+1}'
-\sum_{n=1}^{\infty }2a(-1)^n(\beta_n-1)\beta_{n-1}b_n'b_{n-1}',
\ee
where we have used the relationship $b_{-n}'=(-1)^nb_n'$ given in   (\ref{oned}). 
Moreover, it follows from   (\ref{eq30}) 
that \bbe\nonumber
\lefteqn{ \sum_{n=-\infty}^\infty (-1)^na(\beta_n-1)b_n'[
\beta_{n+1}b'_{n+1}- \beta_{n-1}b'_{n-1}]}\\
 \nonumber&=&\sum_{n=0}^\infty 2a(-1)^n\frac{(\beta_{n+1}-1)\beta_{n+1}b_{n+1}'^2}{\gamma_{n+1}}
-\sum_{n=1}^{\infty }2a(-1)^n(\beta_{n-1}-1)\beta_{n-1}b_{n-1}'^2\gamma_n
\\ \nonumber
&=&
-\sum_{n=1}^\infty 2a(-1)^n\frac{(\beta_{n}-1)\beta_{n}b_{n}'^2}{\gamma_{n}}
+\sum_{n=0}^{\infty }2a(-1)^n(\beta_{n}-1)\beta_{n}b_{n}'^2\gamma_{n+1}
\\
&=&
\sum_{n=1}^\infty 2a(-1)^n(\beta_{n}-1)\beta_{n}b_{n}'^2\kappa_a d_n
+2a(\beta_{0}-1)\beta_{0}b_{0}'^2\gamma_{1},\label{e200}
\bee
where we have used the identity
\be \frac1{\gamma_n}-\gamma_{n+1}=-\kappa_a d_n 
\ee
defined  by   (\ref{gamma}).
Combining  the equations  (\ref{eq40}), (\ref{eq5}) and   (\ref{e200}), we have
\bbe\nonumber
0&=&\left(\sum_{n=1}^\infty+\sum_{n=-1}^{-\infty}+\sum_{n=0}\right) \! \!\!(-1)^n2\kappa_a(\beta_n\!-\!1)\beta_nb_n'^2
\\&&-\sum_{n=1}^\infty (-1)^n4(\beta_{n}-1)\beta_{n}b_{n}'^2\kappa_a \frac{\beta_n}{\beta_n-1}
\nonumber
-2a(\beta_{0}-1)\beta_{0}b_{0}'^2\gamma_{1}
\\
 \nonumber &=&\sum_{n=1}^\infty (-1)^n\left(4\kappa_a-4\kappa_a\frac{\beta_n}{\beta_{n}-1}\right)\beta_n(\beta_n-1)b_n'^2+2\kappa_a\beta_0(\beta_0-1)b_0'^2
\\&&-2a(\beta_{0}-1)\beta_{0}b_{0}'^2\gamma_{1}\nonumber
\\
  &=&-4\kappa_a\sum_{n=1}^\infty (-1)^n\beta_nb_n'^2-2\kappa_a\beta_0b_0'^2,
\label{e300}
\bee
since
\be
-2a(\beta_{0}-1)\beta_{0}b_{0}'^2\gamma_{1}=2a(\beta_{0}-1)\beta_{0}b_{0}'^2\frac{\kappa_a a}{1-a^2}
=-2\kappa_a\beta_0^2b_0'^2
\ee
due to  (\ref{eq40}), (\ref{gamma}) and (\ref{frac}).

On the other hand, multiplying  the $n$th equation of (\ref{eq1}) by $(\beta_n-1)b_n'/(4\kappa_a)$ and summing  the resultant equations yield
\bbe\nonumber
0&=&\frac12\sum_{n=-\infty}^\infty(\beta_{n}-1)\beta_{n}b_{n}'^2
\\
&=&\label{eq10} \sum_{n=1}^\infty(\beta_{n}-1)\beta_{n}b_{n}'^2+\frac12(\beta_0-1)\beta_0 b_0'^2.
\bee
Multiplying   (\ref{e300}) by  $(\beta_0-1)/(4\kappa_a)$ and then adding the resultant equation to    (\ref{eq10}),
we have
\be 0&=&\sum_{n=1}^\infty(\beta_{n}-1)\beta_{n}b_{n}'^2-\sum_{n=1}^\infty (-1)^n(\beta_0-1)\beta_nb_n'^2
\\
&\geq &\sum_{n=2}^\infty(\beta_{n}+\beta_0-2)\beta_{n}b_{n}'^2+(\beta_1-1)\beta_1b_1'^2+(\beta_0-1)\beta_1b_1'^2
\\
&= &\sum_{n=2}^\infty(2a^2+n^2-2)(a^2+n^2)b_{n}'^2+(2a^2-1)(a^2+1)b_1'^2,
\ee
and so, after the use of the condition $2a^2\geq 1$,
\be
0= \sum_{n=2}^\infty(2a^2+n^2-2)(a^2+n^2)b_{n}'^2.
\ee
This implies    $b'_n= 0$  for $n\geq 2$. Substitution of this finding  into (\ref{eq1}) with $n=2$ and $1$
 produces the result $b'_1=0$ and $b_0'=0$.
 Consequently, the validity of the spectral simplicity  expressed by  (\ref{one}) is obtained due to
(\ref{one0}).

The proof of Theorem \ref{spectral} is completed.

\section{Bifurcation analysis}

This section is contributed for the proof of the bifurcation assertion of Theorem \ref{main}.
The following   steady-state bifurcation theorem is crucial to approach the result.
\begin{Theorem}\label{kr} (Krasnoselskii \cite{Kra} and
Nirenberg \cite{Nir}) For a  Banach space $X$, a constant value $\kappa_{\rm{crit}}>0$ and  an open neighborhood $\mathcal D$ of the point $(0,\kappa_{\rm{crit}})$ in the Banach space $ X\times [0,\infty)$, let    $M_\kappa$, $N$  and $F$ be the operators with
$$F(\psi,\kappa)=\psi+\kappa M_\kappa\psi +N(\psi,\kappa),\,\,\, (\psi,\kappa)\in \mathcal{D},$$
subject to the following conditions:
\begin{description}
\item[(i)] $F: \mathcal{D} \mapsto X$ is continuous,
\item[(ii)] $ M_\kappa: \mathcal{D}\mapsto X$ is linear,  compact and continuous,
\item[(iii)] $N: \mathcal {D}\mapsto  X$ is nonlinear  and compact,
\item[(iv)] $N(0,\kappa)\equiv 0$ and $N(\psi,\kappa)=o(\|\psi\|_X)$
uniformly for $(\psi,\kappa)\in \mathcal{D}$,
\item[(v)] the spectral simplicity condition
$$\dim \bigcup_{n=1}^\infty \left\{ \psi\in X, (
Id-\kappa_{\rm{crit}}L_{\kappa_{\rm{crit}}})^n\psi=0\right\}  =1$$
holds true for $Id$ the identity operator in $X$.
\end{description}
 Then there exists a continuous family $(\psi_\kappa,\kappa)\in \mathcal {D}$, different to the trivial one  $(0,\kappa)$, such that
 \bbe\label{e400}F(\psi_\kappa,\kappa)=0,
 \bee
 or the solution family of  (\ref{e400}) branches  off $(0,\kappa_{\rm{crit}})$
when   $\kappa$ varies across the critical value  $\kappa_{\rm{crit}}$.
\end{Theorem}

\

\noindent \textbf{Proof of  Theorem \ref{main}}.
From Lemmas \ref{alpha} and   \ref{cor} we see that a solution $\psi$ bifurcating from $\psi^*$ is regular
whenever $\psi \in
C^{2+\alpha}_{\rm{per}}({\Omega_a})$. Thus it suffices to seek bifurcating solutions in the function  space $\hat C^{2+\alpha}_{\rm{per}}({\Omega_a})$ for $0<\alpha <1$. Recall $\psi^*=\cos x_2$,
 the operator
\be F(\psi,\kappa) =\psi+\kappa
\Delta^{-1}[\kappa+(\nabla\times \psi)\cdot \nabla ]^{-1}\psi^*,
\ee
set in the previous section,
 the operator $L_\kappa$ defined by  (\ref{Lk}) and the critical number  $\kappa_a$  in Theorem \ref{spectral}. For a constant $\epsilon$ such that $ 0<\epsilon <\kappa_a$,  we introduce the
symbols  $$\left\{\begin{array}{l}
X= \hat C^{2+\alpha}_{\rm{per}} ({\Omega_a}),\vspace{2mm}
\\
\mathcal{D}=\left\{\psi\in \hat C^{2+\alpha}_{\rm{per}} ({\Omega_a});\,
\|\nabla^2\psi-\nabla^2\psi^*\|_{C^\alpha} <\frac{(\kappa_a-\epsilon)^2}{4}\right \}\times (\kappa_a-\epsilon,\kappa_a+\epsilon),
\vspace{2mm}
 \\
 M_k\psi = \frac1\kappa(L_k\psi-\psi )=\frac1\kappa \Delta^{-1}[\kappa +(\nabla\times \psi^*)\cdot \nabla  ]^{-1}(\nabla\times \psi^*)\cdot \nabla\psi,
\vspace{2mm} \\
 N(\psi,\kappa)  =F(\psi,\kappa)-\left[\psi-\psi^* +kM_k(\psi-\psi^*)\right].
\end{array}\right.$$
To verify the bifurcation assertion now  remains to demonstrate the validity of
 the assumptions of Krasnolselskii's
theorem.

\emph{Firstly},  we verify   the assumptions (i, ii) of Theorem \ref{kr}.
For the even function property of  $F(\psi,\kappa)(\x)$ with $(\psi,\kappa)\in \mathcal{D}$, we
see that  the even function  $\psi$ implies $\nabla\times\psi$ to be an odd function and so  $\y$. This observation
implies that   $$\psi^*(\y(-\x,s))=\cos(-y_2(\x,s))=\psi^*(\y(\x,s)),$$
and hence  $F(\psi,\kappa)(\x)$ is an even function of $\x \in {\Omega_a}$.

To show the continuity of $F$, for $(\psi,\kappa), \,(\psi',\kappa')\in \mathcal {D}$, we note that
\be
\lefteqn{|\Delta F(\psi',\kappa')-\Delta F(\psi,\kappa)-(\Delta \psi'-\Delta \psi)|}\nonumber\\
 &\leq &  \left|\kappa' [\kappa'+(\nabla\times\psi')\cdot \nabla]^{-1}\psi^* -\kappa'[\kappa+(\nabla\times\psi)\cdot \nabla ]^{-1}\psi^*\right|\nonumber
 \\
 &&+\left|\kappa' [\kappa+(\nabla\times\psi)\cdot \nabla]^{-1}\psi^* -\kappa[\kappa+(\nabla\times\psi)\cdot \nabla ]^{-1}\psi^*\right|\nonumber
  \\
  &= & \kappa' |[\kappa'+(\nabla\times\psi')\cdot \nabla]^{-1}(\kappa-\kappa'+[\nabla\times \psi-\nabla\times \psi']\cdot \nabla) [\kappa+(\nabla\times\psi)\cdot \nabla ]^{-1}\psi^*|\nonumber
  \\
  &&+|\kappa-\kappa'|\,| [\kappa+(\nabla\times\psi)\cdot \nabla]^{-1}\psi^*|.\label{e500}
  \ee
  By the Lagrangian formulation
  \bbe \label{e700} [\kappa+(\nabla\times\psi)\cdot \nabla]^{-1}f(\x)= \int^\infty_0\e^{-\kappa s}f(\y(\x,s))ds,
  \bee
  we have
   \be
\lefteqn{\|\Delta F(\psi',\kappa')-\Delta F(\psi,\kappa)-(\Delta \psi'-\Delta \psi)\|_{C^0}}
  \\
  &\leq&\|(\kappa-\kappa'+[\nabla\times \psi-\nabla\times \psi']\cdot \nabla )[\kappa+(\nabla\times\psi)\cdot \nabla ]^{-1}\psi^*\|_{C^0}+ \frac{|\kappa-\kappa'|}{\kappa}\|\psi^*\|_{C^0}
  \\ &\leq &\frac2{\kappa}|\kappa-\kappa'|+  \|\nabla \times \psi-\nabla\times \psi'\|_{C^0}
  \,\int^\infty_0 e^{-\kappa s} \|\nabla\y\cdot \nabla_\y \psi^*(\y) \|_{C^0}ds,
  \ee
  which is bounded by,  using   (\ref{aaaa2}) and  $\sqrt{\|\nabla^2\psi-\nabla^2\psi^*\|_{C^0}}\leq (\kappa_a-\epsilon)/2$,
  \be
\lefteqn{ \frac{2|\kappa-\kappa'|}{\kappa}+\frac{\left(\sqrt{2}(\kappa-\kappa_a+\epsilon)+\sqrt{5}\right) \|\nabla \psi^*\|_{C^0}\|\nabla\psi-\nabla\psi'
\|_{C^0}}{(\kappa-\kappa_a+\epsilon)^2} }
\\ &\leq& \frac{2|\kappa-\kappa'|}{\kappa}+\frac{\left(\sqrt{2}(\kappa-\kappa_a+\epsilon)+\sqrt{5}\right)\|\nabla\psi-\nabla\psi'
\|_{C^0}}{(\kappa-\kappa_a+\epsilon)^2}.
\end{eqnarray*}
Additionally, by Lemma \ref{alpha}, we have
\begin{eqnarray*}
\|\Delta F(\psi,\kappa)-\Delta \psi\|_{C^1}
&\leq&\kappa \|[\kappa+(\nabla\times\psi)\cdot \nabla]^{-1}\psi^*\|_{C^1}
\\
&\leq &1 +
\frac{\sqrt{2}\kappa(\kappa-2 \sqrt{\|\nabla^2\psi-\nabla^2\psi^*\|_{C^0}})+\sqrt{5}\kappa}{(\kappa-2 \sqrt{\|\nabla^2\psi-\nabla^2\psi^*\|_{C^0}})^2}
\\
&\leq &1+
\frac{\sqrt{2}(\kappa-\kappa_a+\epsilon)\kappa+\sqrt{5}\kappa}{(\kappa-\kappa_a+\epsilon)^2}.
\ee
With the use of the above $C^0$ and $C^1$ estimates, the required continuity of the operator $F$ in the intermediate
H\"older space  is thus derived  from  the interpolation inequality
\be \|\Delta f\|_{C^\alpha}&\leq &2\|\Delta f\|_{C^0}^{1-\alpha} \|\Delta f\|_{C^1}^\alpha
\ee
and the H\"older inequality of the Laplace operator
\be
\|\nabla^2f\|_{C^\alpha}\leq c \|\Delta f\|_{C^\alpha}.
\ee

For  the assumption (ii) of Theorem
\ref{kr}, we rewrite the operator  $M_\kappa$ as
 \be M_k\psi &=&
 \frac1\kappa \Delta^{-1}[\kappa +(\nabla\times \psi^*)\cdot \nabla  ]^{-1}(\nabla\times \psi^*)\cdot \nabla\psi
\\
&=&\frac{1}{\kappa}\Delta^{-1}\psi
-\Delta^{-1}[\kappa+(\nabla\times\psi^*)\cdot\nabla ]^{-1}\psi.
\ee
This formulation enables us to apply the argument on the continuity of the operator $F(\psi,\kappa)$ to obtain the continuity of the operator   $M_\kappa: \mathcal{D}\mapsto X$  and result of
$M_\kappa\psi \in C^{2+\delta}({\Omega_a})$  for any $\alpha<\delta <1$. The  compactness of the operator $M_\kappa$ is due to the compact imbedding of
$C^{2+\delta}({\Omega_a})$ into $C^{2+\alpha}({\Omega_a})$.

\emph{Next}, to verify  the assumptions (iii, iv), we notice that
\be N(\psi,\kappa)
&=& \kappa\Delta^{-1}[\kappa+(\nabla\times\psi)\cdot \nabla ]^{-1}\psi^*+\psi^* -\kappa M_\kappa(\psi-\psi^*).
\ee
 Therefore the compactness of the operator $N$ is implied in the proof of the continuity of $F$ and the compactness of the operator $M_\kappa$.  To prove the non-linear  assertion, we transform the operator $N$ into an explicit quadratic form. That is,  by the solution property of $\psi^*$ satisfying (\ref{bh1}),
\begin{eqnarray*}
N&=& F(\psi,\kappa)-\psi+\psi^*-\kappa M_\kappa(\psi-\psi^*)
\\
&=&
\kappa \Delta^{-1}[\kappa +(\nabla\times\psi)\cdot \nabla]^{-1}\psi^*
\\
&&+\psi^*-\Delta^{-1}[\kappa+(\nabla\times \psi^*)\cdot \nabla
]^{-1}(\nabla\times \psi^*)\cdot \nabla (\psi-\psi^*)
\\
&=&
\kappa\Delta^{-1}[\kappa +(\nabla\times\psi)\cdot \nabla]^{-1}\psi^*-\kappa\Delta^{-1}[\kappa +(\nabla\times\psi^*)\cdot \nabla]^{-1}\psi^*
\\
&&+\Delta^{-1}[\kappa+(\nabla\times \psi^*)\cdot \nabla
]^{-1}(\nabla\times \psi-\nabla\times \psi^*)\cdot \nabla \psi^*,
\end{eqnarray*}
since
$$(\nabla\times \psi^*)\cdot \nabla (\psi-\psi^*)=-(\nabla\times \psi-\nabla\times \psi^*)\cdot \nabla \psi^*.$$
By  elementary manipulations, we have
\begin{eqnarray*}
N&=&
- \kappa\Delta^{-1}[\kappa+(\nabla\times \psi)\cdot \nabla
]^{-1}(\nabla\times \psi-\nabla\times \psi^*)\cdot\nabla (\kappa+(\nabla\times \psi^*)\cdot \nabla)^{-1}\psi^*
\\
&&+\Delta^{-1}[\kappa+(\nabla\times \psi^*)\cdot \nabla
]^{-1}(\nabla\times \psi-\nabla\times \psi^*)\cdot \nabla \psi^*
\\
 &=&-\Delta^{-1}[\kappa+(\nabla\times \psi)\cdot \nabla
]^{-1}(\nabla\times \psi-\nabla\times \psi^*)\cdot\nabla\psi^*
\\
&&+\Delta^{-1}[\kappa+(\nabla\times \psi^*)\cdot \nabla
]^{-1}(\nabla\times \psi-\nabla\times \psi^*)\cdot \nabla \psi^*
\\
 &=&-\Delta^{-1}\left[(\kappa+(\nabla\times \psi)\cdot \nabla
)^{-1}\!-\!(\kappa+(\nabla\times \psi^*)\cdot \nabla
)^{-1}\right](\nabla\times \psi\!-\!\nabla\times \psi^*)\!\cdot\!\nabla\psi^*
\\
&=&\Delta^{-1}[\kappa+(\nabla\times \psi)\cdot \nabla
]^{-1}(\nabla\times \psi-\nabla\times \psi^*)
\\
&&\cdot \nabla [\kappa+(\nabla\times \psi^*)\cdot \nabla
]^{-1}(\nabla\times \psi-\nabla\times \psi^*)\cdot\nabla\psi^*.
\end{eqnarray*}
With the use of this quadratic form and (\ref{e700}), we have
\be
\lefteqn{\|\Delta N(\psi,\kappa)\|_{C^0}}\\
&\leq&\frac1\kappa\|(\nabla\times \psi-\nabla\times \psi^*)\cdot \nabla [\kappa+(\nabla\times \psi^*)\cdot \nabla
]^{-1}(\nabla\times \psi-\nabla\times \psi^*)\cdot\nabla\psi^*\|_{C^0}
\\
&\leq&\frac1\kappa \|\nabla \psi-\nabla\psi^*\|_{C^0}\, \|\nabla [\kappa+(\nabla\times \psi^*)\cdot \nabla
]^{-1}(\nabla\times \psi-\nabla\times \psi^*)\cdot\nabla\psi^*\|_{C^0}
\\
&\leq&\frac1\kappa\|\nabla \psi-\nabla\psi^*\|_{C^0}\, \|\nabla [(\nabla\times \psi-\nabla\times \psi^*)\cdot\nabla\psi^*]\|_{C^0}\int^\infty_0 \e^{-\kappa s}\|\nabla \y^*\|_{C^0}ds
\\
&\leq&\frac1\kappa\|\nabla \psi-\nabla\psi^*\|_{C^0} (\|\nabla^2\psi-\nabla^2\psi^*\|_{C^0}+\|\nabla \psi-\nabla\psi^*\|_{C^0})\int^\infty_0 \e^{-\kappa s}(2+s)ds,
\end{eqnarray*}
where the flow trajectory $\y^*$ is defined by the velocity $\nabla\times \psi^*=(-\sin x_2,0)$ and is in the following form
\bbe\label{e600}
\y^*(\x,t)=\x +t(\sin x_2,\,0).
\bee
 Hence
 $$\|\Delta N(\psi,\kappa)\|_{C^0}\leq
\frac{2\kappa+1}{\kappa^3} \|\psi-\psi^*\|_{C^2}^2.$$

For the estimate of the operator $N$ in the H\"older semi-norm,
we employ   (\ref{aaaa2}), (\ref{e700}) and (\ref{e600}) to produce the estimates
\be
[(\kappa+(\nabla \times \psi)\cdot \nabla
)^{-1}f]_{C^\alpha}&\leq &[f]_{C^\alpha}\int^\infty_0\e^{-\kappa s} \|\nabla \y\|_{C^0}^\alpha ds
\\&\leq &[f]_{C^\alpha}\frac{\sqrt{2}(\kappa-2\alpha\sqrt{\|\nabla^2\psi-\nabla^2\psi^*\|_{C^0}})+\sqrt{5}}{(\kappa-2\alpha\sqrt{\|\nabla^2\psi-\nabla^2\psi^*\|_{C^0}})^2}
\ee
and
\be
\lefteqn{[\nabla (\kappa+(\nabla \times \psi^*)\cdot \nabla
)^{-1}f]_{C^\alpha}}
\\
&=& \int^\infty_0 \e^{-\kappa s}[\nabla \y^*\cdot \nabla _{\y^*} f(\y^*(\cdot,s)]_{C^\alpha}ds
\\
&\leq &\|\nabla f\|_{C^0}\int^\infty_0\e^{-\kappa s} [\nabla \y^*]_{C^\alpha} ds
+[\nabla f]_{C^\alpha}\int^\infty_0\e^{-\kappa s} \|\nabla \y^*\|_{C^0}^{1+\alpha}ds
\\
&\leq&\|\nabla f\|_{C^0}\int^\infty_0\e^{-\kappa s} 2s ds
+[\nabla f]_{C^\alpha}\int^\infty_0\e^{-\kappa s} (2+s)^{1+\alpha}ds
\\
&\leq &\frac2{\kappa^2}\|\nabla f\|_{C^0}
+\frac{4\kappa^2+4\kappa+2}{\kappa^3}[\nabla f]_{C^\alpha}.
\ee
  Let $c$ be a constant independent of $\psi$ and $\kappa$ close to $\kappa_a$ and the constant may change from line to line. Hence for  \be \mathbf{w}=\nabla\times \psi-\nabla\times \psi^*,
 \ee
 the H\"older  semi-norm of the operator $N$ is estimated as
\be
\lefteqn{[ \Delta N(\psi,\kappa)]_{C^\alpha}}\\
&=&[(\kappa+(\nabla\times \psi)\cdot \nabla
)^{-1}\w
\cdot \nabla (\kappa+(\nabla\times \psi^*)\cdot \nabla
)^{-1}\w\cdot\nabla\psi^*]_{C^\alpha}
\\
&\le&c[\w
\cdot \nabla (\kappa+(\nabla\times \psi^*)\cdot \nabla
)^{-1}\w\cdot\nabla\psi^*]_{C^\alpha}
\\
&\le&
c \|\w\|_{C^0}[
 \nabla (\kappa+(\nabla\times \psi^*)\cdot \nabla
)^{-1}\w\cdot\nabla\psi^*]_{C^\alpha}
\\
&&+
c[\w]_{C^\alpha}\|
 \nabla (\kappa+(\nabla\times \psi^*)\cdot \nabla
)^{-1}\w\cdot\nabla\psi^*\|_{C^0}
\\
&\le&
c\left(
\|\w\|_{C^0} \|\nabla(\w\cdot\nabla\psi^*)\|_{C^\alpha}
+[\w]_{C^\alpha}\|\nabla(\w\cdot\nabla\psi^*)\|_{C^0}\right)
\\
&\le&
c
\|\nabla \psi-\nabla \psi^*\|_{C^{1+\alpha}}^2.
\ee
This shows that the assumptions (iii,iv) of Theorem \ref{kr} hold true.

\emph{Finally,} for the verification of the spectral condition, we apply  Theorem
\ref{spectral} to obtain the existence of critical value $\kappa_{a}$ satisfying the simplicity condition
\be\dim\bigcup_{i=1}^\infty\left\{\psi=\sum_{n=-\infty}^\infty b_n\cos (a x_1+nx_2) \in \hat
 C^{2+\alpha}_{\rm{per}}({\Omega_a});\,\,L_{\kappa_{a}}^i \psi=0\right\}= 1.\ee
Therefore, this together with  (\ref{zero})  for $m\neq 1$ produces the validity of the assumption (v) of Theorem \ref{kr}:
\bbe\label{one111}\dim\bigcup_{i=1}^\infty\left\{\psi \in \hat C^{2+\alpha}_{\rm{per}}({\Omega_{a}});\,\,L_{\kappa_a}^i \psi=0\right\}= 1.\bee
The bifurcation assertion of Theorem \ref{main} is thus follows from Theorem \ref{kr} and the proof of  Theorem \ref{main} is completed.

\

\noindent\textbf{Ackownledgement}. The present research  was discussed with  Professor Shouhong Wang when the author  was visiting  the Department of Mathematics at Indiana University  Bloomington in  2002 and the validity of the main result was  basically demonstrated  during the visit. The author  acknowledges  the hospitality  offered by Professor Wang and his group.

\

\

\end{document}